\documentclass[11pt]{amsart}
\usepackage{amsxtra}
\usepackage{amssymb}
\usepackage[mathscr]{eucal}
\usepackage{amssymb,amscd,amsthm,amsmath,graphicx,color}
\usepackage{booktabs}
\usepackage{color}
\usepackage[color,matrix,arrow]{xy}

\usepackage{musicography}

\usepackage{todonotes}
\theoremstyle{plain}

\newtheorem{theorem}{Theorem}[section]
\newtheorem{lemma}[theorem]{Lemma}
\newtheorem{definition-theorem}[theorem]{Definition-Theorem}
\newtheorem{proposition}[theorem]{Proposition}
\newtheorem{corollary}[theorem]{Corollary}
\newtheorem{definition}[theorem]{Definition}
\newtheorem{example}[theorem]{Example}
\newtheorem{remark}[theorem]{Remark}
\newtheorem{conjecture}[theorem]{Conjecture}
\newtheorem{notation}[theorem]{Notation}

\newtheorem{assumption}[theorem]{Assumption}

\newtheorem*{maintheorem*}{Main Theorem}

\newcommand \bth[1] { \begin{theorem}\label{t#1} }
\newcommand \ble[1] { \begin{lemma}\label{l#1} }

\newcommand \bpr[1] { \begin{proposition}\label{p#1} }
\newcommand \bco[1] { \begin{corollary}\label{c#1} }
\newcommand \bde[1] { \begin{definition}\label{d#1}\rm }
\newcommand \bex[1] { \begin{example}\label{e#1}\rm }
\newcommand \bre[1] { \begin{remark}\label{r#1}\rm }
\newcommand \bcj[1] { \begin{conjecture}\label{j#1}\rm }

\newcommand \bnota[1] { \begin{notation}\label{n#1}\rm }
\renewcommand {\eth} { \end{theorem} }
\newcommand {\ele} { \end{lemma} }

\newcommand {\epr} { \end{proposition} }
\newcommand {\eco} { \end{corollary} }
\newcommand {\ede} { \end{definition} }
\newcommand {\eex} { \end{example} }
\newcommand {\ere} { \end{remark} }
\newcommand {\ecj} { \end{conjecture} }

\newcommand {\enota} { \end{notation} }
\newcommand \thref[1]{Theorem \ref{t#1}}
\newcommand \leref[1]{Lemma \ref{l#1}}
\newcommand \prref[1]{Proposition \ref{p#1}}
\newcommand \coref[1]{Corollary \ref{c#1}}

\newcommand \exref[1]{Example \ref{e#1}}

\makeatletter
\newsavebox{\@brx}
\newcommand{\llangle}[1][]{\savebox{\@brx}{\(\m@th{#1\langle}\)}%
  \mathopen{\copy\@brx\kern-0.5\wd\@brx\usebox{\@brx}}}
\newcommand{\rrangle}[1][]{\savebox{\@brx}{\(\m@th{#1\rangle}\)}%
  \mathclose{\copy\@brx\kern-0.5\wd\@brx\usebox{\@brx}}}
\makeatother

\usepackage{tikz}
\usetikzlibrary{matrix}
\usepackage{graphicx,ctable,booktabs}
\usetikzlibrary{arrows.meta}
\usepackage{tikz-cd}

\usepackage{hyperref}
\usepackage{cleveref}

\DeclareMathOperator{\CG}{\mathbf C \rtimes G}

\DeclareMathOperator{\Pic}{Pic} 
 
\DeclareMathOperator{\Spec}{Spec}

 \DeclareMathOperator{\Proj}{Proj}

\DeclareMathOperator{\Aut}{\sf Aut}

 \DeclareMathOperator{\Hom}{Hom}

\DeclareMathOperator{\modd}{{\sf mod}}

\DeclareMathOperator{\coh}{{\sf coh}} 

\DeclareMathOperator{\Spc}{Spc}

\DeclareMathOperator{\Id}{Id}

\DeclareMathOperator{\ev}{ev}
\DeclareMathOperator{\coev}{coev}

\DeclareMathOperator{\SL}{SL}

\DeclareMathOperator{\For}{For}
\DeclareMathOperator{\cone}{cone}
\DeclareMathOperator{\Ind}{Ind}

\DeclareMathOperator{\Perf}{\sf perf}
\DeclareMathOperator{\Cat}{\sf Cat}

\newcommand{\mf}{\mathfrak}
\newcommand{\mc}{\mathcal}

\newcommand{\id}{\operatorname{id}}

\newcommand{\kk}{\Bbbk}

\newcommand{\bR}{\mathbf R}
\newcommand{\bS}{\mathbf S}
\newcommand{\bC}{\mathbf C}
\newcommand{\bK}{\mathbf K}

\newcommand{\bP}{\mathbf P}
\newcommand{\bQ}{\mathbf Q}

\newcommand{\bI}{\mathbf I}
\newcommand{\bJ}{\mathbf J}

\newcommand{\unit}{\ensuremath{\mathbf 1}}

\newcommand{\idem}{\musNatural{}}

\newcommand{\KrG}{\bK \rtimes G}

\newcommand{\GSpc}{G\text{-}\Spc}

\newcommand{\ul}{\underline}

\newcounter{listequation}

\numberwithin{equation}{section}

\begin{document}

\title[Group actions on monoidal triangulated categories]{Group actions on monoidal triangulated categories and Balmer spectra}

\author[Huang]{Hongdi Huang}
\address{(Huang) Department of Mathematics, Rice University, Houston, TX 77005, U.S.A.}
\email{h237huan@rice.edu}

\author[Vashaw]{Kent B. Vashaw}
\address{(Vashaw) Department of Mathematics,
MIT,
Cambridge, MA 02139, U.S.A.}
\email{kentv@mit.edu}



\maketitle

\begin{abstract}
Let $G$ be a group acting on a left or right rigid monoidal triangulated category $\bK$ which has a Noetherian Balmer spectrum. We prove that the Balmer spectrum of the crossed product category of $\bK$ by $G$ is homeomorphic to the space of $G$-prime ideals of $\bK$, give a concrete description of this space, and classify the $G$-invariant thick ideals of $\bK$. Under some additional technical conditions, we prove that the Balmer spectrum of the equivariantization of $\bK$ by $G$ is also homeomorphic to the space of $G$-prime ideals. Examples of stable categories of finite tensor categories and perfect derived categories of coherent sheaves on Noetherian schemes are used to illustrate the theory.

\end{abstract}

\section{Introduction}
Classifications of thick ideals for tensor-triangulated categories have been a problem of major importance for the past 30 years, since the pioneering work of Hopkins \cite{Hopkins1}, Neeman \cite{Neeman2}, Thomason \cite{Thomason1}, Benson--Carlson--Rickard \cite{BCR3}, Friedlander--Pevtsova \cite{FP1}, and others. These problems, arising in the disparate areas of representation theory, algebraic geometry, and homotopy theory, were put into a common framework by Balmer in the early 2000s \cite{Balmer1}: for a tensor-triangulated category $\bK$, there exists a topological space, the Balmer spectrum $\Spc \bK$, whose Thomason-closed sets parametrize the thick ideals of $\bK$. A tensor-triangulated category by definition is a monoidal triangulated category with isomorphisms $A \otimes B \cong B \otimes A$ for all objects $A$ and $B$; for generally noncommutative monoidal triangulated categories, a version of the Balmer spectrum was introduced in \cite{BKS1, NVY1}.

Many computations of Balmer spectra for monoidal triangulated categories have involved some sort of group action on the category in question, and in many examples the Balmer spectrum is realized as a quotient of another Balmer spectrum by a group action, see \cite[Theorem 5.2.2]{BKN1}, \cite[Theorem 9.3.2]{NVY1}, \cite[Theorem 10.4]{NP}, \cite[Theorem 9.1.1]{NVY3}. However, a systematic approach to the situation where a group $G$ acts on a monoidal triangulated category $\bK$ has been lacking; the present paper is an attempt to address this issue directly. A group action on a monoidal category $\bC$ (see Section \ref{group-prelim} for the precise definition) is a choice of monoidal autoequivalences
\[
T_g: \bC \to \bC
\]
for each $g \in G$, together with natural isomorphisms
\[
T_g \circ T_h \cong T_{gh}
\]
satisfying standard compatibility properties. When a group $G$ acts on a finite-dimensional Hopf algebra $H$ by Hopf automorphisms, we get just such an action of $G$ on the tensor category of finite-dimensional modules $\modd(H)$ of $H$. When $\bK$ is monoidal triangulated, we additionally require that each $T_g$ is an equivalence of monoidal triangulated categories. When $H$ is a finite-dimensional Hopf algebra with an action of $G$ as above, then there is an action of $G$ on the monoidal triangulated category $\ul{\modd}(H)$, the stable category of $\modd(H)$. 

The starting point of tensor-triangular geometry is the analogy between a tensor-triangulated category and a ring; the tensor product corresponds to multiplication. If $G$ is a group acting on a ring $R$, then there is a natural way of defining a new ring, called the crossed product of $R$ by $G$, where multiplication is skewed by the action of $G$; closely related are the familiar construction of Ore extensions. Prime ideals of crossed product rings have received significant attention in noncommutative ring theory, see e.g.~ \cite{GM, LP, Passman1, Passman2, Passman3, Goodearl0}.

When a group $G$ acts on a monoidal triangulated category $\bK$, we can analogously define the crossed product category $\bK \rtimes G$ (see Section \ref{crossprod}). Crossed product categories have recently been used by Bergh--Plavnik--Witherspoon to show that every finite tensor category satisfying some mild assumptions can be embedded in a finite tensor category where the tensor product property for cohomological support varieties fails \cite{BPW1}. Examples of crossed product categories include the smash coproduct Hopf algebras studied by Benson--Witherspoon and Plavnik--Witherspoon \cite{BW1, PW}. Balmer spectra for Benson--Witherspoon smash coproducts, as well as some additional Plavnik--Witherspoon smash coproducts, were determined in \cite{NVY1,NVY3}, and it was shown in \cite{Vashaw} that Balmer spectra for certain Benson--Witherspoon smash coproducts coincide with the Balmer spectra of their Drinfeld centers. In all of these cases, the group $G$ acting on the monoidal triangulated category was assumed to be finite. In this paper, we consider the case that $G$ is possibly infinite.

In ring theory, in addition to the crossed product, the invariant subring also plays a critical role in studying group actions on rings (in particular in relation to prime ideals, see e.g.~ \cite{Montgomery, LMS, Taylor}). The analogue in the categorical setting is that of an equivariantization (see Section \ref{sect-equiv}). Objects of the equivariantization are defined to be pairs consisting of an object from the original category, and a family of isomorphisms which trivialize the $G$-action in a compatible way. 

When a group $G$ acts on a ring $R$, the third related concept to the crossed product and invariant subring one considers is that of a $G$-prime ideal, which is defined to be a $G$-invariant ideal which satisfies the prime condition on the level of $G$-invariant ideals. We make the analogous definition in the monoidal triangulated setting, and denote by $\GSpc \bK$ the $G$-Balmer spectrum of $\bK$, that is, the collection of $G$-prime ideals under a suitable topology. The relationship between the prime spectrum of $R$ and the $G$-prime spectrum of $R$ has also attracted significant attention, including the celebrated Goodearl--Letzter Stratification Theorem \cite{GL, Lorenz2}. Motivated by classical ring theory, the following are now natural questions that we wish to explore, when $G$ is a group acting on a monoidal triangulated category $\bK$:
\begin{enumerate}
    \item What is the Balmer spectrum of the crossed product category $\bK \rtimes G$? Does it parametrize the thick ideals of $\bK \rtimes G$?
    \item What is the Balmer spectrum of the equivariantization $\bK^G$? Does it parametrize the thick ideals of $\bK^G$?
    \item What is the $G$-Balmer spectrum $\GSpc \bK$ of $\bK$? Does it parametrize the $G$-invariant thick ideals of $\bK$?
\end{enumerate}

Here we say that a topological space $X$ parametrizes thick ideals when there is a bijection between thick ideals and Thomason-closed or specialization-closed subsets of $X$. {\it{Specialization-closed}} means arbitrary union of closed sets, and {\it{Thomason-closed}} means union of closed sets with quasicompact complement. It is known that the Balmer spectrum parametrizes thick ideals in many cases: when the tensor product of $\bK$ is commutative, more generally when all prime ideals of $\bK$ are completely prime, when $\bK$ has a thick generator, and when $\Spc \bK$ is Noetherian (see \thref{balmer-class} below). Part of the difficulty in question (1) above is that the crossed product category $\KrG$ does not satisfy, a priori, any of these conditions (when $G$ is infinite). One difficulty in question (2) is that in general we do not know a systematic way to endow $\bK^G$ the structure of a triangulated category; however, the triangulated structure exists under a few additional conditions, thanks to work of Elagin \cite{Elagin}. We are able to answer the above questions in the following way in terms of $\Spc \bK$.

\bth{maintheorem}(See \prref{up-down-ideals}, \prref{equivariantization}, \prref{topol-gspc}, \prref{ideal-class-gspc}).
Let $G$ be a group acting on a monoidal triangulated category $\bK$. Then 
\begin{enumerate}
    \item $\Spc (\bK \rtimes G)$ is homeomorphic to $\GSpc \bK$;
    \item if $G$ is finite, $\bK$ is $\kk$-linear with $|G| \not = 0$ in $\kk$, and $\bK^G$ is canonically triangulated (e.g. $\bK$ is the stable category of a finite tensor category, or more generally admits a DG enhancement), then $\Spc \bK^G$ is homeomorphic to $\GSpc \bK$;
    \item if $\bK$ is left or right rigid and $\Spc \bK$ is Noetherian, then specialization-closed subsets of $\GSpc \bK$ are in bijection with $G$-invariant thick ideals, points of $\GSpc \bK$ are in bijection with $G$-orbits of points of $\Spc \bK$, and closed sets of $\GSpc \bK$ are generated by all orbits of closed sets of $\Spc \bK$.
\end{enumerate}
\eth

The assumptions on $\bK$ in part (3) of this theorem are not too restrictive. Indeed, in our primary motivating example, $\bK=\ul{\modd}(H)$ for a finite-dimensional Hopf algebra $H$, we have that $\bK$ is rigid and it is a conjecture that $\Spc \bK$ is Noetherian (see \exref{stablecat} below for details).

\subsection*{Acknowledgements}
The authors thank Pavel Etingof, Ken Goodearl, Dan Nakano, Cris Negron, and Sarah Witherspoon for helpful conversations and comments. We also thank Greg Stevenson for extensive discussions, and in particular for pointing out the argument in \exref{perf-action}. Some results in this paper were formulated at the Structured Quartet Research Ensembles (SQuaREs) program in March 2022 and April 2023 at the American Institute of Mathematics, and during a visit by the second author to the Department of Mathematics at Rice University in March 2023; we thank the AIM, and the second author thanks Rice University, for their hospitality and support. The second author was partially supported by the NSF postdoctoral fellowship DMS-2103272.

\section{Monoidal triangular preliminaries}\label{sec:s1}

We review the basics of monoidal triangulated categories and their underlying geometry. A {\em{monoidal triangulated category}} is a triangulated category (c.f.~ \cite{Neeman}) with a compatible tensor product $\otimes$ and unit $\unit$, as in \cite{NVY1}. Throughout, we will denote the shift on a triangulated category $\bK$ by $\Sigma: \bK \to \bK$. We do not require any braiding condition on the monoidal product. We have standard definitions below which follow the noncommutative generalization \cite{NVY1} of Balmer's tensor triangular geometry \cite{Balmer1}.
\begin{enumerate}
\item A subcategory $\bI$ of $\bK$ is called a {\em{triangulated subcategory}} if for any distinguished triangle
    \[
    A \to B \to C \to \Sigma A,
    \]
    if $A$ and $B$ are in $\bI$, then so is $C$, and $A$ is in $\bI$ if and only if $\Sigma A$ is in $\bI$.
\item A triangulated subcategory $\bI$ is called {\em{thick}} if $A \oplus B \in \bI$ implies $A$ and $B$ are both in $\bI$.
\item A thick subcategory $\bI$ is called a {\em{two-sided thick ideal}} if $A \in \bI$ implies $A \otimes B$ and $B \otimes A$ are in $\bI$, for any object $B \in \bK$. ``Two-sided thick ideal" will often be shortened to ``thick ideal." 
\item Given a collection of objects $\bS$, we write $\langle \bS \rangle$ for the smallest thick ideal containing $\bS$.
\item A thick ideal $\bI$ of $\bK$ is called {\em{principal}} if it can be written $\bI=\langle A \rangle$ for some object $A$. If an ideal has a finite generating set, then it is principal, since the smallest thick ideal containing a direct sum of objects is the smallest thick ideal containing the union of those objects.
\item A proper thick ideal $\bP$ will be called {\em{prime}} if $A \otimes \bK \otimes B \subseteq \bP$ implies $A$ or $B \in \bP$, or equivalently if $\bI \otimes \bJ \subseteq \bP$ implies that either $\bI$ or $\bJ \subseteq \bP$, for all thick ideals $\bI$ and $\bJ$ of $\bP$ \cite[Theorem 3.2.2]{NVY1}.
\item {\em{The Balmer spectrum of $\bK$}}, denoted $\Spc \bK$, is the topological space consisting of the set of prime thick ideals, where closed sets have the form
    \[
    V_{\bK}(\bS) := \{\bP \in \Spc \bK : \bS \cap \bP = \varnothing\},
    \]
    where $\bS$ is a collection of objects of $\bK$. When the category $\bK$ is clear from context, we write $V_{\bK}(\bS):=V(\bS).$
\item The {\it{support}} of an ideal $\bI$ in $\bK$ is defined as 
\[
\Phi_V (\bI):=\bigcup_{A \in \bI} V(A)=\{\bP \in \Spc \bK : \bI \not \subseteq \bP\}.
\] By definition, this is a specialization-closed subset of $\Spc \bK$.
\item We call $\bK$ {\it{left rigid}} if every object  $A$ of $\bK$ has a left dual $A^*$, as in \cite[Section 2.10]{EGNO}. {\it{Right rigid}} monoidal triangulated categories are defined similarly. We call $\bK$ {\it{half-rigid}} if it is either left or right rigid.
\item A thick ideal is called {\it{semiprime}} if it can be written as an intersection of prime ideals. When $\bK$ is half-rigid, all thick ideals are semiprime \cite[Proposition 4.1.1]{NVY2}.
\item We call $\bK$ {\it{$\Spc$-Noetherian}} if $\Spc \bK$ is Noetherian, as a topological space. 
\end{enumerate}

In many cases, there is a bijection between Thomason-closed subsets of $\Spc \bK$ and thick ideals of $\bK$. We recall these results now for reference.
\bth{balmer-class}
Let $\bK$ be a monoidal triangulated category. In each of the following cases, the support map $\Phi_V$ induces an order-preserving bijection between semiprime thick ideals of $\bK$ and Thomason-closed subsets of $\Spc \bK$:
\begin{enumerate}
    \item when $A \otimes B \cong B \otimes A$ for all objects $A$ and $B$ of $\bK$ \cite[Theorem 4.10]{Balmer1};
    \item when all prime ideals of $\bK$ are completely prime (that is, when $A \otimes B \in \bP$ implies $A$ or $B$ is in $\bP$, for any objects $A$ and $B \in \bK$ and any prime ideal $\bP$ of $\bK$) \cite[Theorem 3.11]{MallickRay};
    \item when $\Spc \bK$ is Noetherian \cite[Theorem 4.9]{Rowe};
    \item when $\bK$ is has a thick generator \cite[Theorem A.7.1]{NVY4}.
\end{enumerate}
When $\bK$ is half-rigid, the bijection in each case restricts to a bijection between principal ideals and closed subsets of the form $V(A)$, for $A \in \bK$.
\eth

There are also various lattice-theoretic instantiations of \thref{balmer-class}, see \cite[Proposition]{BKS1}, \cite[Corollary 10]{Krause}, \cite[Theorem 6.4.5]{GS}. The bijection between semiprime ideals and specialization-closed sets can be used to show the following.

\bco{noeth-desc}
Let $\bK$ be a half-rigid $\Spc$-Noetherian monoidal triangulated category. Then principal ideals of $\bK$ satisfy the descending chain condition.
\eco

We end the section by reminding the reader of some of the standard monoidal triangulated categories that we are interested in studying, which will serve as running examples throughout the paper.

\bex{stablecat}
Let $\bC$ be a finite tensor category, that is, $\bC$ is an abelian rigid monoidal category with finite-dimensional morphism spaces, finitely many simple objects, enough projectives, and the unit is simple. For example, $\bC=\modd(H)$, the category of finite-dimensional modules for a finite-dimensional Hopf algebra $H$, is a finite tensor category. The stable category of $\bC$, denoted $\ul{\bC}$, is defined as having the same objects as $\bC$, and morphism spaces are quotients by those morphisms factoring through projective objects. Then $\ul{\bC}$ is a rigid monoidal triangulated category (c.f.~ \cite{Happel1}). There is a conjectural description of the Balmer spectrum of $\ul{\bC}$ \cite[Conjecture E]{NVY3}. In particular, $\Spc \ul{\bC}$ is conjectured to be Noetherian.
\eex

\bex{perf}
Let $X$ be a topologically Noetherian scheme. Denote by $\Perf(X)$ the derived category of perfect complexes of coherent sheaves on $X$ (c.f.~ \cite[Section 3.2.3]{Rouquier2010}). Note that $\Perf(X)$ is the subcategory of compact objects in the full unbounded derived category of $\coh(X)$, \cite[Lemma 3.5]{Rouquier2010} or \cite[Theorem 3.1.1]{BV}. Then $\Perf(X)$ is a rigid symmetric monoidal triangulated category, and its Balmer spectrum is homeomorphic to $X$, by \cite[Theorem 6.3(a)]{Balmer1}, using work of Thomason \cite{Thomason1}. The Thomason isomorphism will be denoted 
\[
\eta_X: X \xrightarrow{\cong} \Spc \Perf(X).
\]
When $R$ is a commutative ring, we will typically abbreviate $\Perf(R):=\Perf(\Spec R),$ and the Thomason isomorphism as $\eta_R: \Spec R \to \Spc \Perf(R)$. 
\eex

\section{Group action preliminaries}
\label{group-prelim}

We now recall the definition of a group action on a monoidal category, following closely \cite[Definition 4.13]{DGNO} and \cite[Section 2.7]{EGNO}. Let $\bC$ be a monoidal category and $G$ a group. We will typically denote $e$ for the identity element of $G$, and will denote the group operation in $G$ using multiplicative notation. An {\it action of $G$ on $\bC$} consists of the following data:
\begin{enumerate}
    \item for each $g \in G$, a monoidal autoequivalence $T_g: \bC \to \bC$;
    \item for each pair $g, h \in G$, a natural isomorphism $\gamma_{g,h} : T_g T_h \xrightarrow{\cong} T_{gh}$, such that for any three elements $g,h, k \in G$, the diagram
\begin{center}
    \begin{tikzcd}
T_{g}T_{h} T_{k} \arrow[d, "\gamma_{gh} T_k"'] \arrow[r, "T_g \gamma_{hk}"] & T_g T_{hk} \arrow[d, "{\gamma_{g,hk}}"] \\
T_{gh} T_k \arrow[r, "{\gamma_{gh,k}}"']                                    & T_{ghk}                                
\end{tikzcd}
\end{center}
    commutes.
\end{enumerate}
\bre{te-id}
Note that it is a consequence of (1) and (2) that $T_e$ is naturally isomorphic to the identity functor $\Id:\bC \to \bC$ (\cite[Remark 3.2]{Elagin}).
\ere

These conditions can alternatively be packaged together in the following way. To the group $G$, we define a monoidal category $\Cat(G)$ where objects are elements of $G$, the only morphisms are identity morphisms $g \to g$, and the monoidal product is given by multiplication in $G$. Recall we also have a monoidal category $\Aut_{\otimes}(\bC)$ where objects are monoidal autoequivalences of $\bC$, and the monoidal product is the composition of functors. Then an action of $G$ on $\bC$ is just a monoidal functor 
\[
\Cat(G) \to \Aut_{\otimes}(\bC).
\]

Now let $\bK$ be a monoidal triangulated category. In this case, we require each $T_g$ to be a triangulated functor; in other words, we replace the category of monoidal autoequivalences $\Aut_{\otimes}(\bC)$ with the category of autoequivalences of monoidal triangulated categories. Recall that this means we have natural isomorphisms $T_g \Sigma \cong \Sigma T_g$ for each $g \in G$ such that if
\[
A \to B \to C \to \Sigma A
\]
is a distinguished triangle, then the corresponding triangle
\[
T_g A \to T_g B \to T_g C \to \Sigma T_g A
\]
is also distinguished.

We now give a few examples of actions of groups on monoidal triangulated categories that are of natural interest.

\bex{fintens-action}
Suppose $G$ acts on a finite tensor category $\bC$ (c.f.~ \exref{stablecat}). Since any equivalence of abelian categories is exact and preserves projectives, it is straightforward that there is an induced $G$-action on $\ul{\bC}$. 
\eex

\bex{hopfalg-action}
In particular, if $H$ is a finite-dimensional Hopf algebra, and $G$ a group acting on $H$ by Hopf algebra automorphisms, then $G$ acts on $\modd(H)$ via setting $T_g(M)$ to be the module with $H$-action $*$ given in terms of the $H$-action $.$ on $M$ by $h*m:=g^{-1}(h).m$. In particular, if $G$ is an algebraic group and $H$ is an $\mc{O}(G)$-comodule Hopf algebra, then $G$ acts on $H$ by Hopf automorphisms. A natural example is the Drinfeld double $D(\mc{O}(G)^*)$ of the dual of the coordinate ring of a finite group scheme $G$, which is a $\mc{O}(G)$-comodule Hopf algebra. 
\eex

\bex{cocleft-action}
Not all group actions on $\modd(H)$ arise from actions at the level of algebras, however; if $H$ is a cocleft extension of $\mc{O}(G)$, then $G$ also acts on $\modd(H)$. For example, by \cite[Section 2.3]{DEN}, if $G$ is a simple algebraic group over an algebraically closed field $\kk$ of characteristic 0 with Lie algebra $\mf{g}$, then $G$ acts on the category $\modd(u_{\zeta}(\mf{g}))$; this action does not arise from an action on the level of Hopf algebras. Here $u_{\zeta}(\mf{g})$ is Lusztig's small quantum group associated to a primitive $\ell$th root of unity $\zeta$, where $\ell$ is an odd integer compatible with the Cartan data of $G$ \cite{Lusztig}. 
\eex

\bex{perf-action}
Let $G$ be a discrete group acting on a topologically Noetherian scheme $X$, and consider the action of $G$ on $\Perf(X)$ via derived pullback, as in \cite[Example 3.4]{Elagin}. The induced action of $G$ on the $\Spc \Perf(X)$ corresponds to the original action of $G$ on $X$ via the Thomason isomorphism $\eta_X: X \to \Spc \Perf(X)$, by an argument of Stevenson \cite{Stevenson} which we include here.

Let $x$ be a point of $X$. Pick an affine open $\Spec \mc{O}_X(U) \cong U \subseteq X$ which contains $x$. We have an equivalence of monoidal triangulated categories
\[
(\Perf(X)/\bI(X\backslash U))^{\idem} \cong \Perf(U),
\]
where $\bI(X \backslash U)$ is the ideal consisting of objects in $\Perf(X)$ supported on the complement of $U$, by \cite[Theorem 2.13]{Thomason1} (see also \cite[Lemmas 3.4 and 3.10]{Rouquier2010}), induced by the restriction functor $R_U: \Perf(X) \to \Perf(U)$. Here $(-)^{\idem}$ is the idempotent-completion operation (see \cite{Balmer-Schlichting}); note however that, working on the level of $\Spc$, the idempotent completion does not affect one's calculation \cite[Corollary 3.14]{Balmer1}. It follows that  
\begin{align}
    \label{etaxu}
    \eta_X(x) & = \{A \in \bK: R_U A \in \eta_U(x)\}.
\end{align}

Recall that by \cite[Proposition 8.1]{Balmer2}, we have a functorial comparison map
\begin{align*}
    \Spc \Perf(U) &\xrightarrow{\rho_U} U \cong \Spec \mc{O}_X(U) \\
    \bP &\mapsto \langle f \in \mc{O}_X(U) : \cone(f) \not \in \bP\rangle,
\end{align*}
which is an isomorphism (indeed, the inverse of the Thomason isomorphism $\eta_U$). Here we use implicitly the identification 
\[
\mc{O}_X(U) \cong \Hom_{\Perf(U)} (\mc{O}_U, \mc{O}_U)
\]
when we refer to $\cone(f)$ for $f \in \mc{O}_X(U)$. We also have a functor $R_x: \Perf(U) \to \Perf( \mc{O}_{X,x})$ induced by the ring homomorphism $\mc{O}_X(U)\to \mc{O}_{X,x}$. By the functoriality of $\Spc$ for braided monoidal triangulated categories \cite[Proposition 3.6]{Balmer1}, we have induced maps on spectra
\[
\Spc \Perf(\mc{O}_{X,x}) \to \Spc \Perf(U) \to \Spc \Perf(X).
\]
By the functoriality of the comparison map, we have a commutative diagram
\begin{center}
    \begin{tikzcd}
{\Spc \Perf(\mc{O}_{X,x})} \arrow[d] \arrow[r, "\rho_{\mc{O}_{X,x}}"] & {\Spec \mc{O}_{X,x}} \arrow[d] \\
\Spc \Perf(U) \arrow[r, "\rho_{U}"]                        & U,                             
\end{tikzcd}
\end{center}
where the vertical arrows are the usual maps. The right hand side of this diagram sends the unique closed point of $\Spec \mc{O}_{X,x}$ to $x \in U$. The top row of the diagram sends the 0-ideal of $\Perf(\mc{O}_{X,x})$ to the unique closed point in $\Spec \mc{O}_{X,x}$, since $\Perf(X)$ is rigid, and hence the unique closed point of $\Spc \Perf(\mc{O}_{X,x})$ must be the 0-ideal. The left column, by definition, sends the 0-ideal to the collection of objects $\bP:=\{A \in \Perf(U) : R_x A=0\}$. By commutativity of the diagram, we have $\rho_U(\bP)=x$, that is, $\bP=\eta_U(x)$. It follows from (\ref{etaxu}) that 
\[
\eta_X(x) = \{ A \in \Perf(X) : R_U A \in \bP \} = \{A \in \Perf(X) : R_x'A \cong 0 \},
\]
where we denote $R_x' = R_x R_U: \Perf(X) \to \Perf( \mc{O}_{X,x})$ the standard restriction functor. Let $g \in G$. Now we can compute
\begin{align*}
T_g\eta_X(x) &= T_g \bP'\\
&= \{ T_g A : R_x' A \cong 0 \}\\
&=\{A: R_x' T_{g^{-1}}A \cong 0\}\\
&= \{A : R_{g.x}'A \cong 0\}\\
&= \eta_X(g.x),
\end{align*}
that is, the Thomason isomorphism is $G$-equivariant. Here the fourth equality follows from the commutative diagram 
\begin{center}
    \begin{tikzcd}
\Perf(X) \arrow[r,"T_g"] \arrow[d, "R_x'"]    & \Perf(X) \arrow[d, "R_{g.x}'"]    \\
{\Perf(\mc{O}_{X,x})} \arrow[r] & {\Perf(\mc{O}_{X,g.x})}
\end{tikzcd}
\end{center}
where the bottom row is the functor induced by the action of $G$ on $X$. 
\eex

\bex{picard-action}
Let $\bK$ be any rigid monoidal triangulated category. Recall the Picard group $\Pic \bK$ of $\bK$, defined as the collection of isomorphism classes tensor-invertible objects of $\bK$ \cite[Definition 2.1]{Balmer2.5}, with group operation given by the tensor product. Then $\Pic \bK$ acts on $\bK$ by conjugation. That is, if $A$ is a tensor-invertible object of $\bK$, i.e.~ $A \otimes A^* \cong \unit$, then denoting the isomorphism class of $A$ in $\Pic \bK$ by $[A]$, the action of $\Pic \bK$ on $\bK$ is given by $
T_{[A]} B := A \otimes B \otimes A^*.$
\eex

\section{$G$-ideals and $G$-prime ideals}

Throughout this section, let $G$ be a group acting on a monoidal triangulated category $\bK$. We introduce some of the primary objects of study in this paper, namely $G$-ideals and $G$-primes, and prove some elementary properties which will be useful going forward. 

It is straightforward to verify that if $\bI$ is a thick ideal of $\bK$ and $g \in G$, then 
\[
T_g \bI:=\{T_g A : A \in \bI\}
\]
is another thick ideal, since $T_g$ is a monoidal triangulated functor. 

\bde{g-prime}
A thick ideal $\bI$ of $\bK$ will be called a {\em{$G$-ideal}} if it is $G$-invariant, that is, $T_g \bI \subseteq \bI$ for all $g \in G$. A $G$-ideal $\bQ$ will be called a {\em{$G$-prime}} if for any $G$-invariant ideals $\bI$ and $\bJ$ of $\bK$, we have $\bI \otimes \bJ \subseteq \bQ$ implies $\bI$ or $\bJ \subseteq \bQ$. The collection of $G$-primes of $\bK$ will be denoted $\GSpc \bK$, and be considered as a topological space where closed sets are 
\[
V^G(\bS) := \{ \bQ \in \GSpc \bK : \bS \cap \bQ = \varnothing\}
\]
for any collection $\bS$ of objects in $\bK$. We have the $G$-Balmer support $\Phi_G$, which sends a $G$-ideal $\bI$ to the specialization-closed subset
\[
\Phi_G (\bI):= \bigcup_{A \in \bI} V^G(A)
\]
of $\GSpc \bK$.
\ede

\bre{cont-equ}
Since thick ideals are closed under isomorphism, and since there is a natural isomorphism $T_g T_{g^{-1}} \cong T_e \cong \Id$, if $\bI$ is a $G$-ideal then $T_g \bI = \bI$ for all $g \in G$.
\ere

\bex{picard-action-ideal}
Let $\bK$ be a rigid monoidal triangulated category, and consider the action of $G:=\Pic \bK$ on $\bK$ as in \exref{picard-action}. For $[A] \in \Pic \bK$, we have $T_{[A]} B \in \langle B \rangle$, and it follows that every thick ideal of $\bK$ is a $G$-ideal of $\bK$, hence $\Spc \bK \cong \GSpc \bK$. 
\eex

The following lemma gives a natural way of constructing $G$-ideals of $\bK$. 

\ble{g-ideal-gens}
Let $G$ be a group acting on a monoidal triangulated category $\bK$. Let $\bS$ be a $G$-invariant collection of objects of $\bK$. Then $\langle \bS \rangle$ is a $G$-ideal of $\bK$.
\ele

\begin{proof}
Set 
\[
\bI:=\{ A \in \langle \bS \rangle : T_g A \in \langle \bS \rangle \;\; \forall g \in G \}.
\]
Since $\bS$ is itself $G$-invariant, it follows that $\bS \subseteq \bI$. We claim that $\bI$ in fact is a thick ideal, and hence $\langle \bS \rangle \subseteq \bI$; since the opposite containment holds by definition, it will follow that $\bI = \langle \bS \rangle$. It is clear that $\bI$ is a thick subcategory, using the fact that each $T_g$ is a triangulated functor; for example, if 
\[
A \to B \to C \to \Sigma A
\]
is a distinguished triangle with $A$ and $B$ both in $\bI$, then since 
\[
T_g A \to T_g B \to T_g C \to \Sigma T_g A
\]
is a distinguished triangle with both $T_gA$ and $T_gB$ in $\langle \bS \rangle$, it follows that $T_g C \in \langle \bS \rangle$, that is, $C \in \bI$. The ideal property for $\bI$ follows from the fact that each $T_g$ is monoidal: if $A \in \bI$ and $B \in \bK$, then we see that 
\[
T_g(A \otimes B) \cong T_g(A) \otimes T_g(B)
\]
is in $\langle \bS \rangle$ since $T_gA \in \langle \bS \rangle$, and it follows that $A \otimes B \in \bI$. 
\end{proof}

Note that if $\bP$ is a prime ideal of $\bK$, then $T_g \bP$ is also a prime ideal of $\bK$. Hence we have a natural action of $G$ on $\Spc \bK$, via
\[
g.\bP :=T_g \bP = \{T_g A : A \in \bP\}.
\]
We now connect $G$-invariance on the level of ideals with $G$-invariance on the level of $\Spc.$

\ble{g-ideal-supp}
Let $G$ be a group acting on a monoidal triangulated category $\bK$. If $\bI$ is a $G$-ideal, then the support $\Phi_V(\bI)$ of $\bI$ is a $G$-invariant subset of $\Spc \bK$. If $\bI$ is a semiprime ideal, and if $\Phi_V(\bI)$ is $G$-invariant, then $\bI$ is a $G$-ideal. 
\ele

\begin{proof}
Recall that by \Cref{sec:s1} (8)
$$\Phi_V(\bI) = \{ \bP \in \Spc \bK : \bI \not \subseteq \bP\}.$$
If $\bI$ is $G$-invariant, then $\bI \not \subseteq \bP$ for $\bP$ a prime ideal implies that $\bI \not \subseteq T_g\bP$ for any $g \in G$, that is, if $\bP$ is in $\Phi_V(\bI)$ then $T_g\bP$ is in $\Phi_V(\bI)$ as well. Thus $\Phi_V(\bI)$ is $G$-invariant.

On the other hand, if $\bI$ is semiprime and $\Phi_V(\bI)$ is $G$-invariant, and if $\bP \supseteq \bI$, then $T_g\bP\supseteq \bI$ for all $g\in G$ as well. Hence,
since
\[
\bI = \bigcap_{\bP' \supseteq \bI} \bP',
\]
if $A \in \bI$ then $A \in \bP'$ for all $\bP'$ over $\bI$, hence by the above $T_g A$ is in $\bP'$ for all $\bP'$ over $\bI$, and so $T_g \bI \subseteq\bI$ for all $g\in G$. 
\end{proof}

We give a useful reformulation of the $G$-prime condition. 

\ble{gprime-char}
Let $G$ be a group acting on a monoidal triangulated category $\bK$. Let $\bQ$ be a $G$-ideal of $\bK$. Then the following are equivalent:
\begin{enumerate}
\item $\bQ$ is a $G$-prime;
\item $\bQ$ satisfies
\[
T_g A \otimes \bK \otimes T_h B \subseteq \bQ \; \forall \; g,h \in G \Rightarrow A \text{ or }B \in \bQ;
\]
\item $\bQ$ satisfies
\[
T_g A \otimes \bK \otimes  B \subseteq \bQ \; \forall \; g \in G \Rightarrow A \text{ or }B \in \bQ;
\]
\item $\bQ$ satisfies
\[
A \otimes \bK \otimes T_h B \subseteq \bQ \; \forall \; h \in G \Rightarrow A \text{ or }B \in \bQ.
\]
\end{enumerate}
\ele

\begin{proof}
We first check that (1) and (2) are equivalent. Suppose $\bQ$ is $G$-prime, and that $T_g A \otimes \bK \otimes T_h B \subseteq \bQ$ for some objects $A$ and $B \in \bK$, and all $g$ and $h \in G$. Then it is straightforward that
\[
\langle T_g A : g \in G \rangle \otimes \langle T_h B : h \in G \rangle \subseteq \bQ.
\]
Since both ideals are $G$-ideals by \leref{g-ideal-gens}, it follows that either $A$ or $B$ is in $\bQ$. 

For the other direction, assume $\bQ$ satisfies (2), and suppose that 
\[
\bI \otimes \bJ \subseteq \bQ
\]
for two $G$-ideals $\bI$ and $\bJ$. Assume $\bI \not \subseteq \bQ$; then there exists some $A \in \bI$ with $A \not \in \bQ$. But then
\[
T_g A \otimes \bK \otimes T_h B \subseteq \bI \otimes \bJ \subseteq \bQ
\]
for all $B \in \bJ$ and $g$ and $h \in G$, and by assumption either $A$ or $B$ is in $\bQ$. Since $A \not \in \bQ$ by assumption, $B \in \bQ$ for all $B \in \bJ$, that is, $\bJ \subseteq \bQ$. It follows that (1) and (2) are equivalent.

It is clear that (3) implies (2), by definition. To see that (2) implies (3), just note that if 
\[
T_g A \otimes \bK \otimes B \subseteq \bQ
\]
for all $g$ in $G$, then using the fact that $\bQ$ is a $G$-ideal,
\[
T_{hg} A \otimes \bK \otimes T_h B \subseteq \bQ
\]
for all $g, h \in G$, and so by (2) either $A$ or $B$ is in $\bQ$. The equivalence of (2) and (4) is similar.
\end{proof}

Finally, we prove that under a Noetherian assumption, a principal or prime ideal is invariant under a functor $T_g$, for $g \in G$, if and only if it is equal to its image under that functor. A ring theory analogue of the lemma below may be found in \cite[Section 2]{Goodearl0}. 

\ble{G-inv}
Let $G$ be a group acting on a half-rigid $\Spc$-Noetherian monoidal triangulated category $\bK$. If $\bI$ is a principal ideal so that $T_g \bI \subseteq \bI$ for some $g \in G$, then $T_g \bI= \bI$.
\ele
\begin{proof}
Since $\bI$ is principal, $T_g \bI$ is principal. Thus we obtain a descending chain of principal ideals
\[
\bI \supseteq T_g \bI \supseteq T_{g^2} \bI \supseteq ... \supseteq T_{g^n} \bI ...
\]
Since principal ideals of $\bK$ satisfy DCC by \coref{noeth-desc}, there exists $n$ with $T_{g^n}(\bI) = T_{g^{n+1}} (\bI)$. By applying the equivalence $T_{g^{-n}}$ to this identity, the result is obtained. 
\end{proof}

\ble{tg-prime}
Let $G$ be a group acting on a half-rigid $\Spc$-Noetherian monoidal triangulated category $\bK$. If $\bP \in \Spc \bK$ satisfies $T_g \bP \subseteq \bP$ for some $g \in G$, then $T_g \bP = \bP$. 
\ele

\begin{proof}
    Let $\bP$ be a prime with $T_g \bP \subseteq \bP$. Since $\Spc \bK$ is Noetherian and $\bK$ is half-rigid, by \thref{balmer-class} there exists an object $A$ with 
    \[
    V(A):=\{\bP'\in \Spc\bK : A \not \in \bP'\} = \overline{\{\bP\}} = \{\bP'\in \Spc\bK : \bP' \subseteq \bP\}.
    \]
    Since $T_g \bP \subseteq \bP$, we have $A \not \in T_g\bP$, in other words, $T_{g^{-1}}A \not \in \bP$. That is, $\bP \in V(T_{g^{-1}}A)$. Hence \[
    \overline{\{\bP\}} = V(A) \subseteq V(T_{g^{-1}} A).
    \]
    Now again by \thref{balmer-class}, since there is a containment preserving bijection between principal ideals of $\bK$ and closed subsets of the form $V(B)$, we have $A \in \langle T_{g^{-1}} A \rangle.$ By \leref{G-inv}, we have $\langle A \rangle = \langle T_{g^{-1}} A \rangle$, and so
    \begin{align*}
        \overline{\{\bP\}} &= V(A) \\
        &= V(T_{g^{-1}} A)\\
        &= \overline{\{T_{g^{-1}} \bP\}}.
    \end{align*}
    This implies that $\bP = T_{g^{-1}} \bP$, and by applying the autoequivalence $T_g$ we observe that that $\bP = T_g \bP$. 
\end{proof}

\section{Ideals and primes in the crossed product category}
\label{crossprod}

Let $\bC$ be a monoidal category, and $G$ be a group acting on $\bC$. As in \cite[Section 4.15]{EGNO}, we set the \emph{crossed product category} of $\bC$ by $G$ to be a direct sum
\[
\bC \rtimes G := \bigoplus_{g \in G} \bC
\]
indexed by elements of $G$, as additive categories. If $A$ is an object of $\bC$, then the corresponding object of $\bC \rtimes G$ in the copy of $\bC$ indexed by $g$ will be denoted by $A \boxtimes g$ (the notation is meant to reflect the fact that if $\bC$ is a $\kk$-linear category, then $\bC \rtimes G$ is given by the familiar Deligne tensor product \cite[Section 1.11]{EGNO}). That is, objects of $\bC \rtimes G$ are formal direct sums of objects of the form $A \boxtimes g$, for $A \in \bC$ and $g \in G$, and $\Hom_{\bK \rtimes G} (A \boxtimes g, B \boxtimes h):=\Hom_{\bK}(A,B)$ if $g=h$, and 0 otherwise. A monoidal product is then defined via
\[
(A \boxtimes g) \otimes (B \boxtimes h):=(A \otimes T_g (B)) \boxtimes gh.
\]
The monoidal unit for $\bC \rtimes G$ is the object $\unit \boxtimes e$, where $e$ is the identity of $G$.

\ble{duals-cross}
Let $G$ be a group acting on a monoidal category $\bC$. If the object $A$ in $\bC$ has a left (resp.~ right) dual, then $A \boxtimes g$ has a left (resp.~ right) dual in $\bC \rtimes G$. In particular, if $\bC$ is half-rigid, then so is $\CG$.
\ele

\begin{proof}
We give the argument for left duals. Let $A$ be in $\bC$ with left dual $A^*$. We claim that $A \boxtimes g$ has left dual $T_{g^{-1}} (A^*) \boxtimes g^{-1}$. In particular, we define evaluation by
\[
(T_{g^{-1}} (A^*) \boxtimes g^{-1}) \otimes (A \boxtimes g) \cong T_{g^{-1}}(A^* \otimes A) \boxtimes e \xrightarrow{T_{g^{-1}}( \ev_A) \boxtimes e} \unit \boxtimes e
\]
and coevaluation by
\[
\unit \boxtimes e \xrightarrow{\coev_A \boxtimes e} (A \otimes A^*) \boxtimes e \cong  (A \boxtimes g) \otimes (T_{g^{-1}} (A^*) \boxtimes g^{-1}).
\]
Now note that
\[
A \boxtimes g \xrightarrow{ (\coev_A \boxtimes e)  \otimes \id_{A \boxtimes g}} (A \boxtimes g) \otimes (T_{g^{-1}} (A^*) \boxtimes g^{-1}) \otimes (A \boxtimes g) \xrightarrow{\id_{A \boxtimes g} \otimes T_{g^{-1}} (\ev_A) \boxtimes e}  A \boxtimes g
\]
corresponds to the morphism
\[
A \boxtimes g \xrightarrow{(\coev_A \otimes \id_A) \boxtimes g} (A \otimes A^* \otimes A) \boxtimes g \xrightarrow{(\id_A \otimes \ev_A) \boxtimes g} A \boxtimes g
\] by the definition of tensor product of morphisms in $\bC \rtimes G$, and this morphism is $\id_{A \boxtimes g}$ by the duality between $A^*$ and $A$ in $\bC$. This gives one of two identities required for duality; the other follows similarly.
\end{proof}

If $\bK$ is a monoidal triangulated category, and $G$ is a group acting on $\bK$, then we observe that $\KrG$ has the natural structure of a monoidal triangulated category, since it is (as an additive category) simply a direct sum of triangulated categories. 

\bex{fin-tens-cprod}
Let $G$ be a group acting on a finite tensor category $\bC$. As in \exref{stablecat}, we can form the stable category $\ul{\bC}$, which is monoidal triangulated, and the action of $G$ on $\bC$ defines an action of $G$ on $\ul{\bC}$, so we may form $\ul{\bC} \rtimes G$. On the other hand, we may also form $\ul{\bC \rtimes G}$. Since $\bC \rtimes G$ is simply a direct sum of copies of $\bC$ as an additive category, it is straightforward to check that the indecomposable projective objects of $\bC \rtimes G$ are precisely objects of the form $P \boxtimes g$ for $P$ a projective object of $\bC$ and $g \in G$. It follows that the natural functor $\ul{\bC \rtimes G} \to \ul{\bC} \rtimes G$ is an equivalence of monoidal triangulated categories.
\eex

\bex{group-hopf}
Let $H$ be a finite-dimensional Hopf algebra, and $G$ a finite group acting on $H$ by Hopf algebra automorphisms. Recall that $G$ acts on $\modd(H)$ and $\ul{\modd}(H)$ as in \exref{hopfalg-action}. At the abelian level, we have $\modd(H) \rtimes G \cong \modd( (H^* \# \kk G)^*)$. Furthermore, 
\[
\ul{\modd}((H^* \# \kk G)^*) \cong \ul{\modd}(H) \rtimes G
\]
as above in \exref{fin-tens-cprod}.
Here $(H^* \# \kk G)^*$ is the Plavnik--Witherspoon Hopf algebra considered in \cite{BW1,PW}. There is a conjectured description of the Balmer spectrum for the stable category of such a Hopf algebra \cite[Conjecture E]{NVY3}, which has been verified for several large families of examples (\cite[Theorem 9.2.3]{NVY1}, \cite[Theorem 9.1.1]{NVY3}).
\eex

Given an action of a group $G$ on a monoidal triangulated category $\bK$, we now connect $G$-ideals and $G$-primes with ideals and primes in the crossed-product category. We set some notation: if $\bI$ is a $G$-ideal of $\bK$, we denote $\bI \rtimes G$ to be the collection of objects of $\KrG$ generated additively by objects of the form $A \boxtimes g$ for all $A \in \bI$ and $g \in G$. We embed $\bK$ into $\KrG$ by sending $A \mapsto A \boxtimes e$; hence, $\bJ \cap \bK$ will refer to the subcategory of $\bK$ consisting of objects $A$ such that $A \boxtimes e \in \bJ$, if $\bJ$ is a thick ideal of $\KrG$. 

We now check that these maps give bijections between appropriate collections of thick and prime ideals. Parts (1)--(3) may be viewed as an analogue of \cite[Lemma 1.1]{LP}, although the analogue of part (3) for rings is not a bijection, but rather, just a surjection in one direction.

\bpr{up-down-ideals}
Let $G$ be a group acting on a monoidal triangulated category $\bK$. Let $\bI$ be a $G$-ideal of $\bK$ and $\bJ$ be a thick ideal of $\bK \rtimes G$. We have:
\begin{enumerate}
\item $\bI \rtimes G$ is a thick ideal of $\KrG$;
\item $\bJ \cap \bK$ is a $G$-ideal of $\bK$;
\item the maps given by (1) and (2) define a bijection between $G$-ideals of $\bK$ and thick ideals of $\KrG$, that is, $(\mathbf I\rtimes G)\cap \bK=\mathbf I$
and $(\bJ\cap \bK)\rtimes G = \bJ$;
\item $\bQ$ is a $G$-prime of $\bK$ if and only if $\bQ \rtimes G$ is a prime ideal of $\KrG$;
\item the bijection from (3) restricts to a homeomorphism $\Spc (\KrG) \cong \GSpc \bK.$
\end{enumerate}
\epr

\begin{proof}
First we show (1). By the definition of the triangulated structure on $\KrG$, it is straightforward that $\bI\rtimes G$ is triangulated, and in fact thick. To check the ideal property, note that by the definition of the tensor product, 
\[
(A \boxtimes g) \otimes (B \boxtimes h) = (A \otimes T_g(B)) \boxtimes gh
\]
for all $A, B$ in $\bK$ and $g, h \in G$; if $A \boxtimes g \in \bI\rtimes G$, i.e.~ $A \in \bI$, then $A \otimes T_g(B) \in \bI$ and consequently $(A \otimes T_g(B)) \boxtimes gh \in \bI\rtimes G$. On the other hand, if $B \boxtimes h \in \bI\rtimes G$, i.e.~ $B \in \bI$, then $T_g B \in \bI$ since it is a $G$-ideal, and again $(A \boxtimes g) \otimes (B \boxtimes h) = (A \otimes T_g(B)) \boxtimes gh \in \bI\rtimes G$.

For (2), it is similarly straightforward to show that $\bJ \cap \bK$ is a thick ideal of $\bK$. It remains to show that it is $G$-invariant. Let $A \in \bJ \cap \bK$, i.e.~ $A \boxtimes e \in \bJ$. Then 
\[
(\unit \boxtimes g) \otimes (A \boxtimes e) \otimes (\unit \boxtimes g^{-1}) = T_g(A) \boxtimes e \in \bJ,
\]
hence $T_g A \in \bJ \cap \bK$. 

For (3), just note
\begin{align*}
    (\mathbf I\rtimes G)\cap \bK &= \{ \text{direct sums of elements of the form }A \boxtimes g \text{ for }A \in \bI, g \in G\} \cap \bK\\
    &=\{ A : A \boxtimes e  \in \{ \text{direct sums of objects of the form }A \boxtimes g \text{ for }A \in \bI, g \in G\} \}\\
    &= \bI,
\end{align*}
and
\begin{align*}
    (\bJ\cap \bK)\rtimes G &= \{A : A \boxtimes e \in \bJ\} \rtimes G\\
    &= \{\text{direct sums of objects of the form } A \boxtimes g \text{ for } A \text{ with } A \boxtimes e \in \bJ\}\\
    &=\bJ.
\end{align*}
The last equality holds due to the fact that if $\bigoplus_{g\in G}A_g\boxtimes g\in \bJ$, then since $\bJ$ is thick, each $A_g\boxtimes g\in \bJ$ and so $A_g\boxtimes e=(A_g\boxtimes g) \otimes (\mathbf 1\boxtimes g^{-1})\in \bJ$. 

(4) follows from (3) since for any $G$-ideals $\bI$, $\bJ$, and $\bQ$ of $\bK$, we have
\begin{align*}
    \bI \otimes \bJ \subseteq \bQ &\Leftrightarrow (\bI \otimes \bJ) \rtimes G \subseteq \bQ \rtimes G\\
    & \Leftrightarrow (\bI \rtimes G) \otimes (\bJ \rtimes G) \subseteq \bQ \rtimes G.\\
    \end{align*}

For part (5), the only thing we need to verify is that the restriction of the maps from (3) are continuous. Let $\bQ \in \GSpc \bK$ and $\bS$ be a collection of objects of $\KrG$. Then we check
\begin{align*}
    \bQ \rtimes G \in V_{\KrG}(\bS) & \Leftrightarrow \forall \bigoplus_{i \in I} A_i \boxtimes g_i \in \bS, \exists i \text{ with } A_i \boxtimes g_i \not \in \bQ \rtimes G\\
    &\Leftrightarrow \forall  \bigoplus_{i \in I} A_i \boxtimes g_i \in \bS, \exists i \text{ with } A_i \not \in \bQ\\
    & \Leftrightarrow \bQ \in V^G_{\bK}\left ( \left \{\bigoplus_{i \in I} A_i : \bigoplus_{i \in I} A_i \boxtimes g_i \in \bS \text{ for some collection of }g_i \in G \right \} \right).
\end{align*}
On the other hand, if $\bP \in \Spc( \KrG)$ and $\bR$ is a collection of objects in $\bK$, then we similarly check
\begin{align*}
    \bP \cap \bK \in V^G_{\bK}(\bR) & \Leftrightarrow \forall A \in \bR, A \boxtimes e \not \in \bP\\
    &\Leftrightarrow \bP \in V_{\KrG}(A \boxtimes e: A \in \bR).
\end{align*}
\end{proof}

\bco{duals}
Let $G$ be a group acting on a half-rigid monoidal triangulated category $\bK$. Then every $G$-ideal of $\bK$ is an intersection of $G$-primes.
\eco

\begin{proof}
By \leref{duals-cross}, $\KrG$ is half-rigid. Recall that this implies every thick ideal of $\KrG$ is semiprime, i.e.~ is an intersection of prime ideals (see \Cref{sec:s1}\,(10)). The corollary then follows by \prref{up-down-ideals}.
\end{proof}

We end this section with a consequence with regards to the tensor product property on $\GSpc \bK$. The Balmer support $V^G$ is said to have the tensor product property if 
\[
V^G(A \otimes B)= V^G(A) \cap V^G(B)
\]
for all objects $A$ and $B$, or, equivalently, if every prime ideal is completely prime. The tensor product property for various support varieties has been intensely studied recently, see \cite{BW1,PW,NVY2,NP0,BPW0,BPW1}. 

\bco{tpp-action}
Let $G$ be a group acting on a half-rigid $\Spc$-Noetherian monoidal triangulated category $\bK$. If the induced action of $G$ on $\Spc \bK$ is nontrivial on closed points, then $\bK \rtimes G$ has a nonzero nilpotent element. As a consequence, $V^G$ does not satisfy the tensor product property. On the other hand, if all closed points of $\Spc \bK$ are fixed under the action of $G$, then $\bK \rtimes G$ has nilpotent objects if and only if $\bK$ does.
\eco

\begin{proof}
    Suppose $\bP$ is a closed point of $\Spc \bK$ (that is, a minimal prime) such that $T_g \bP \not = \bP$. We have an object $A$ such that $V(A)=\{\bP\}$ by \thref{balmer-class}. Then 
    \[
    V(A \otimes T_g A) \subseteq V(A) \cap V(T_gA) = \{\bP \} \cap \{T_g \bP\} = \varnothing,
    \]
    hence $A \otimes T_gA \cong 0$. Therefore, $A \boxtimes g$ in $\bK \rtimes G$ is nilpotent (of order 2), and by \cite[Theorem 4.2.1]{NVY2}, the Balmer support for $\bK \rtimes G$ does not satisfy the tensor product property. By \prref{up-down-ideals}, $V^G$ does not satisfy the tensor product property. 

    It is clear that if $\bK$ has a nilpotent object, then so does $\bK \rtimes G$. Suppose that $\bK \rtimes G$ has a nonzero nilpotent object, and that the action of $G$ on $\Spc \bK$ fixes all closed points. In particular, then, there exists some $A \boxtimes g$ which squares to 0, hence $A \otimes T_gA \cong 0$. Let $\bP$ be some closed point of $V(A)$. Then note
    \[
    \bigcup_{B \in \bK} V(T_g A \otimes B \otimes A)= V(T_g A) \cap V(A) \supseteq \{\bP\},
    \]
    so there exists some $B$ in $\bK$ so that $C:=T_g A \otimes B \otimes A$ is nonzero. But since $A \otimes T_g A \cong 0$, we have $C\otimes C\cong 0$, so $\bK$ has a nonzero nilpotent.

\end{proof}

\section{Ideals and primes in the equivariantization}
\label{sect-equiv}
Let $\bK$ be a monoidal triangulated category and $G$ a group acting on $\bK$. In this section we connect (under some technical conditions) the $G$-Balmer spectrum of $\bK$ to the Balmer spectrum of $\bK^G$, the equivariantization of $\bK$, which plays the role analogous to the invariant subring in ring theory. 

We begin by recalling the standard definitions, see \cite[Section 4.1.3]{DGNO} or \cite[Section 2.7]{EGNO}. Suppose that $G$ is a group acting on a monoidal category $\bC$ (we will not use any triangulated structure until later in the section). A {\em $G$-equivariant object in $\bC$} is defined as a pair $(A,\phi)$ consisting of $A \in \bC$ and a family of isomorphisms $\phi=\{\phi_g: g\in G\}$ where $\phi_g$ is an isomorphism $T_g A \xrightarrow{\cong} A$, such that the diagrams 
\[
 \xymatrix{
 T_g T_h A\ar[rr]^-{\gamma_{g,h}^A}\ar[d]^-{T_g \phi_h} & &T_{gh} A\ar[d]^-{\phi_{gh}}\\
 T_{g}A\ar[rr]^-{\phi_{g}}&& A,
 }   
\] commute, where $\gamma_{g,h}$ is the natural isomorphism $T_g\circ T_h\cong T_{gh}$. The {\em equivariantization $\bC^G$ of $\bC$} the category where objects are equivariant objects of $\bC$ and morphisms $(A,\phi) \to (B,\psi)$ are defined as morphisms $A \to B$ in $\bC$ which commute with the equivariant structure. The equivariantization inherits a monoidal product from $\bC$ in the straightforward way.

Suppose $G$ is finite. We have two natural functors between $\bC$ and $\bC^G$:
\begin{center}
\begin{tikzcd}
A \arrow[rr, maps to]             &  & {\left ( \bigoplus_{g \in G} T_g A, \tau^A \right )} \\
\bC \arrow[rr, "\Ind", bend left] &  & \bC^G \arrow[ll, "\For", bend left]                  \\
A                                 &  & {(A,\phi)} \arrow[ll, maps to]                      
\end{tikzcd}
\end{center}
 where $\tau^A=\{\tau_h^A : h \in G\}$ is the straightforward equivariant structure given explicitly by 
\[
\tau_h^A = \bigoplus_{g\in G}\gamma_{h, h^{-1}g}^A : T_h \left (\bigoplus_g T_g A \right ) \xrightarrow{\cong}  \bigoplus_g T_g A. 
\]
The functors $\Ind$ and $\For$ are both left and right adjoint to one another (see \cite[Lemma 4.6(a)]{DGNO}). 
We obtain the following familiar tensor identity.
\ble{keyequation}
Let $\bC$ be a monoidal category and $G$ a finite group acting on $\bC$. Suppose $A$ and $(B, \phi)$ are objects in $\bC$ and $\bC^G$, respectively. Then 
\[
\Ind(A \otimes \For(B, \phi))\cong \Ind(A) \otimes (B, \phi) .
\]
\ele
\begin{proof}
We claim that the map $\bigoplus_{g\in G}(\id_{T_gA}\otimes \phi_g)$ from $\Ind(A\otimes \For(B, \phi)) \cong \Ind(A \otimes B)$ to $\Ind(A)\otimes(B, \phi)$ is an isomorphism in $\bC^G$. Indeed, this map is compatible with the equivariant structures, since the diagram
\[
\xymatrix{
T_h \left ( \bigoplus_{g \in G} T_{h^{-1}g} (A \otimes B) \right ) \ar[rrrr]^-{\tau_h^{A \otimes B}}\ar[d]^-{\cong}                   & &  & &\bigoplus_{g \in G} T_g(A \otimes B) \ar[d]^-{\cong} \\
\bigoplus_{g \in G} \left(T_h T_{h^{-1}g}(A)  \otimes T_h T_{h^{-1}g} (B)\right) \ar[rrrr]^-{\bigoplus_{g\in G}\left(\gamma_{h, h^{-1}g}^A \otimes \gamma_{h,h^{-1}g}^B\right)}\ar[d]^-{ \bigoplus_{g\in G} \left(\id_{T_{h}T_{h^{-1}g}(A)} \otimes T_h(\phi_{h^{-1}g})\right)} & &  & & \bigoplus_{g \in G}\left(T_g(A) \otimes T_g(B)\right) \ar[d]^-{\bigoplus_{g\in G} \id_{T_gA} \otimes \phi_g} \\
\bigoplus_{g \in G}  T_h T_{h^{-1} g} (A) \otimes T_h(B)  \ar[rrrr]^-{ \;\;\;\;\;\;\;\;\;\;\;\;\tau_h^A \otimes \phi_h}                              & & & &  \bigoplus_{g \in G} T_g(A) \otimes B  }
\]
commutes: the top square commutes since each $\gamma$ is a natural transformation of monoidal functors (\cite[Definition 2.4.8]{EGNO}), and the bottom square commutes by the fact that $\phi$ is an equivariant structure for $B$. 
\end{proof}

At this point, we will now assume that $\bK$ is monoidal triangulated and $G$ is a finite group acting on $\bK$. We need to make some additional assumptions, in order to ensure that $\bK^G$ has a standard triangulated structure. That is, we wish to call a triangle in $\bK^G$ distinguished if and only if its image under $\For$ in $\bK$ is distinguished. By a theorem of Elagin \cite[Theorem 6.5, Theorem 6.9]{Elagin}, using work of Balmer \cite{Balmer3}, the following assumptions are enough to verify that $\bK^G$ endowed with the distinguished triangles as above is indeed triangulated:
\begin{assumption}
    \label{trikg}
    The group $G$ is finite, the triangulated category $\bK$ is $\kk$-linear for some ring $\kk$ with $|G| \not = 0$ in $\kk$, and $\bK$ is the homotopy category of some stable model category (or, in particular, $\bK$ admits a DG enhancement). 
\end{assumption}

For the remainder of the section, the conditions of Assumption \ref{trikg} are in effect.

\bex{equiv-stable}
Let $G$ be a finite group acting on a finite tensor category $\bC$ over a field $\kk$ as in \exref{stablecat} with $|G|\not = 0$ in $\kk$. The equivariantization $\bC^G$ of $\bC$ by $G$ is again a finite tensor category \cite[Section 4.15]{EGNO}, and indeed $\ul{\bC^G} \cong \ul{\bC}^G$ (here the former category is the stable category of the equivariantization of $\bC$, whereas the latter is the equivariantization of the stable category of $\bC$) \cite[Theorem 3.14]{Sun}. That is, taking the stable category commutes with taking the equivariantization. Note that this situation satisfies Assumption \ref{trikg} since $\ul{\bC}$ has a DG enhancement \cite[Theorem 3.8]{Keller}.
\eex

\bex{equiv-deriv}
Let $G$ be a group acting on a Noetherian scheme $X$ with the induced action on $\Perf(X)$, as in \exref{perf-action}. Assume this action satisfies Assumption \ref{trikg} (note that $\Perf(X)$ does have a DG enhancement, by e.g.~ \cite[Section 3.1]{LS}). Denote by $\bC := \coh(X)$, the category of coherent sheaves on $X$. Then by \cite[Example 3.20]{Sun} (or alternatively by \cite[Remark 7.2]{Elagin} in conjunction with \cite[Proposition 3.2]{BN}), we have $\Perf(X)^G \cong {\sf D}(\bC^G)^c$ the compact part of the derived category of the equivariantization of $\bC$. Under some conditions (see e.g. \cite[Introduction]{Elagin}), the category $\bC^G$ is the category of coherent sheaves on the quotient variety $X/G$; more generally, it can be constructed as the category of coherent sheaves on the quotient stack $X//G$. When $\coh(X)^G \cong \coh(X/G)$, it follows that $\Perf(X)^G \cong \Perf(X/G)$, that is, taking the perfect derived category commutes with taking the equivariantization.
\eex

Suppose $\bI$ is a thick ideal of $\bK$. Then we write 
\[
\bI^G := \{ (A,\phi) \in \bK^G : A \in \bI\}.
\]
On the other hand, suppose $\bJ$ is a thick ideal of $\bK^G$. Then we write
\[
\bJ \cap \bK := \{ A \in \bK : \Ind(A) \in \bJ\}.
\]
Here the notation is meant to reflect the fact that $\Ind$ embeds $\bK$ as a subcategory of $\bK^G$. Via these maps, we obtain a bijection between $G$-ideals of $\bK$ and ideals of $\bK^G$ in the following proposition. In the analogous ring theoretic setting no such bijection exists, but rather there is a bijection between $G$-primes and a certain quotient of the spectrum of the invariant subring \cite[Proposition 4.2]{Montgomery}. 

\bpr{equivariantization}
Let $G$ be a group acting on a monoidal triangulated category $\bK$ satisfying Assumption \ref{trikg} and let $\bI$ be a $G$-ideal of $\bK$ and $\bJ$ be an ideal of $\bK^G$. We have:
\begin{enumerate}
    \item $\bI^G$ is a thick ideal of $\bK^G$;
    \item $\bJ \cap \bK$ is a $G$-ideal of $\bK$;
    \item the maps given by (1) and (2) define a bijection between $G$-ideals of $\bK$ and thick ideals of $\bK^G$, that is, $(\bI^G) \cap \bK = \bI$ and $(\bJ \cap \bK)^G = \bJ;$
    \item $\bP$ is a $G$-prime of $\bK$ if and only if $\bP^G$ is a prime ideal of $\bK^G$;
    \item the bijection from (3) restricts to a homeomorphism $\Spc \bK^G \cong \GSpc \bK.$
\end{enumerate}
\epr
\begin{proof}
(1) is straightforward, since $\For$ is a monoidal triangulated functor. For (2), we have that $\Ind$ is a triangulated functor, but no longer monoidal, so that it is still straightforward to check that $\bJ \cap \bK$ is a thick subcategory; it remains to show the ideal property. Let $A \in \bJ \cap \bK$, and let $B \in \bK$. Since $B$ is a summand of $\For(\Ind(B))$, we have $\Ind(A \otimes B)$ is a summand of $\Ind(A \otimes \For(\Ind(B))) \cong \Ind(A) \otimes \Ind(B)$, using \leref{keyequation}. Since $\Ind(A) \in \bJ$ by assumption, and $\bJ$ is a thick ideal, $\Ind(A) \otimes \Ind(B) \in \bJ$, and hence $\Ind(A \otimes B) \in \bJ$, and it follows that $\bJ \cap \bK$ is a thick ideal. Since $\Ind(A) \cong \Ind(T_g A)$ for all $g \in G$, it follows that $\bJ \cap \bK$ is a $G$-ideal.

(3): We may compute
\begin{align*}
   \bI^G \cap \bK&= \{(A, \phi)\in \bK^G\,:\, A\in \bI\} \cap \bK\\
  &=\{B\in \bK\,:\, \Ind(B)\in\{(A, \phi) \in \bK^G : A \in \bI\}\} \\
  &=\{B \in \bK : \For \Ind (B) \in \bI\}\\
  &=\bI.
\end{align*}
The last step follows from the fact that $\bI$ is a $G$-ideal. For the other direction, note that 
\begin{align*}
    (\bJ \cap \bK)^G&=\{ A \in \bK : \Ind(A) \in \bJ \}^G\\
    &=\{(B, \phi) \in \bK^G : B \in \{ A \in \bK : \Ind(A) \in \bJ \}\}\\
    &=\{(A, \phi) \in \bK^G\,:\, \Ind(A)\in \bJ\};
\end{align*}
we must show that this is equal to $\bJ$. Let $(A, \phi)$ be an object with $\Ind(A)\in \bJ$. We have that $(A,\phi)$ is a summand of $\Ind(A)$ by \cite[Lemma 4.6(b)]{DGNO};  hence $(A, \phi)\in \bJ$, and so 
\[
(\bJ \cap \bK)^G \subseteq \bJ.
\]
Now let $(A,\phi)\in \bJ.$ Then 
\[
\Ind(\unit)\otimes (A, \phi)\cong \Ind(\unit\otimes A)=\Ind(A)
\]
by \leref{keyequation}, and so $\Ind(A) \in \bJ$. It follows that $(\bJ \cap \bK)^G = \bJ$.

Using the bijection from (3), one can check (4) and (5) similarly to \prref{up-down-ideals}(4) and (5). 

\end{proof}

\bre{surj}
When $\bK$ is the stable category of a finite tensor category $\bC$ over a field $\kk$, and $G$ is a finite group with $|G|\not = 0 \in \kk$ acting on $\bC$ as in \exref{fintens-action}, an alternate proof of the surjectivity of the map in \prref{equivariantization}(5) can be approached by the arguments in \cite{BCHS}, using \exref{equiv-stable} and the fact that $\bK$ is the compact part of a compactly-generated monoidal triangulated category (the ``big" stable category, c.f.~ \cite[Appendix A]{NVY3}).
\ere

\section{The core map and a description of $\GSpc$}

Throughout the section, we assume $G$ is a group acting on a monoidal triangulated category $\bK$. In the previous sections, we have connected $\GSpc \bK$ with $\Spc \bK \rtimes G$, and with $\Spc \bK^G$ (under certain assumptions). In this section, we give a complete description of $\GSpc \bK$, in terms of the induced action of $G$ on $\Spc \bK$. We first give analogues of $G$-cores of ideals and prime ideals from ring theory, c.f.~ \cite[Section 11.7.1]{Lorenz1}.

\bde{g-core}
Let $\bI$ be a thick ideal of $\bK$. Then the {\em{$G$-core of $\bI$}}, denoted $(\bI:G)$, is defined as the largest $G$-ideal contained in $\bI$. 
\ede

It is straightforward to verify that
\[
(\bI: G)=\bigcap_{g \in G} T_g \bI.
\]
Here, recall that $T_g \bI$ refers to the thick ideal consisting of elements of the form $T_g A$, for $A \in \bI$. It is also straightforward to observe that if $\bP$ is a prime ideal of $\bK$, then $(\bP: G)$ is a $G$-prime ideal of $\bK$. Hence we have an induced core map
\begin{align*}
    \Spc \bK & \xrightarrow{c} \GSpc \bK\\
    \bP & \mapsto (\bP:G).
\end{align*}
Here we consider $c$ only as a set map.

Just as in the ring-theoretic case (see \cite[Section 0.2]{Lorenz5}), it is clear that $c$ factors through the topological quotient $(\Spc \bK)/G$, whose points are $G$-orbits of primes in $\Spc \bK$, 
via
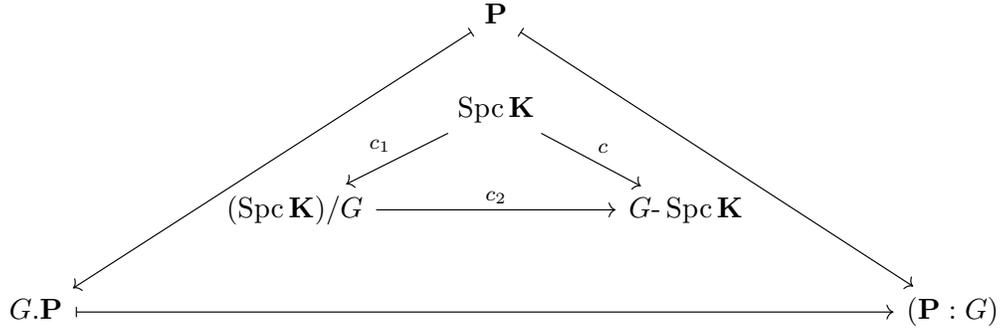
\begin{center}
    \begin{tikzcd}
 &  &                                & \bP \arrow[lllddd, maps to] \arrow[rrrddd, maps to] &           &  &         \\
     &  &   & \Spc \bK \arrow[ld, "c_1"'] \arrow[rd, "c"]         
     &    &  &   \\
 &  & (\Spc \bK)/G \arrow[rr, "c_2"] &                    & \GSpc \bK &  &     \\
G.\bP \arrow[rrrrrr, maps to] &  &   &   &    &  & (\bP:G)
\end{tikzcd}
\end{center}
where $G.\bP$ denotes the orbit of $\bP$ in $\Spc \bK$ under the induced $G$-action.


We next prove the bijectivity of the map $c_2$ when $\Spc \bK$ is Noetherian. This result shows a deviation of tensor-triangular geometry from ring theory, since for even very well-behaved group actions on rings, the map from $G$-orbits of primes to the $G$-prime spectrum typically fails to be injective. On the other hand, surjectivity of the core map is satisfied in ring theory under Noetherianity assumptions, see \cite[p. 338]{GM} and \cite[Lemma II.1.10]{BG} (although the proofs are quite different than the ones we provide in the monoidal triangulated case). In the ring theory case, this gives a bijection between {\em{finite}} orbits of prime ideals and $G$-prime ideals, whereas in the monoidal triangulated case below we obtain a bijection between {\em{arbitrary}} orbits of prime ideals and $G$-primes.

\bpr{bijective-core}
Let $G$ be a group acting on a half-rigid $\Spc$-Noetherian monoidal triangulated category $\bK$. Then the map $c_2: (\Spc \bK)/G \to \GSpc \bK$ given by $G.\bP \mapsto (\bP: G)$ is bijective.
\epr

\begin{proof}
We first show injectivity of $c_2$. For the sake of contradiction, suppose that $c(\bP)=c(\bP')$ for some $\bP, \bP' \in \Spc \bK$ with $G.\bP \not = G.\bP'$. Suppose that $T_g \bP \subseteq \bP'$ for some $g \in G$ and $T_h \bP' \subseteq \bP$ for some $h \in G.$ In that case, we would have that 
\[
T_{hg} \bP \subseteq T_h \bP' \subseteq \bP,
\]
and so by \leref{tg-prime} this would imply 
\[
T_{hg} \bP = T_h \bP' = \bP,
\]
hence that $G.\bP = G.\bP'$, a contradiction. Thus, without loss of generality we may assume that $T_g \bP \not \subseteq \bP'$ for all $g \in G$, that is, for each $g \in G$, we can find $A_g \in \bP$ such that $T_g A_g \not \in \bP'$.

We now construct a descending chain of closed subsets of $\Spc \bK$, by choosing $i$ such that at each step the containment is strict (if it is possible to do so):
\[
V(T_{g_1} A_1) \supsetneq V(T_{g_1}A_1) \cap V(T_{g_2}A_2) \supsetneq ... \supsetneq V(T_{g_1}A_1) \cap V(T_{g_2}A_2) \cap...\cap V(T_{g_i} A_i) \supsetneq...
\]
By the Noetherianity assumption, this process must terminate, say at step $n$. Also by Noetherianity, there exists an object $A$ such that 
\[
V(A) = V(T_{g_1}A_1) \cap V(T_{g_2}A_2) \cap...\cap V(T_{g_n} A_n).
\]
We claim that this $A$ satisfies (1) $A \in c(\bP)$, and (2) $A \not \in c(\bP')$, showing that $c(\bP) \not = c(\bP').$

For (1), we note that by assumption, for any $g\in G$, we have
\[
V(A) \cap V(T_g A_g)=V(A).
\]
In other words, $V(A) \subseteq V(T_g A_g)$, which implies, by the order-preserving bijection between closed subsets and principal ideals \thref{balmer-class}, that $A \in \langle T_g A_g \rangle$. Since $T_g A_g \in T_g \bP$, we have $A \in T_g \bP$ for all $g$, and it follows $A \in c(\bP).$

For (2), to show that $A\not \in c(\bP')$ it suffices to note that $A \not \in \bP'$: if $A$ was in $\bP'$, we would have $\bP' \not \in V(A).$ But that would imply that
\[
\bP' \not \in V(T_{g_1}A_1) \cap V(T_{g_2}A_2) \cap...\cap V(T_{g_n} A_n),
\]
which would mean by definition that for some $i = 1 \cdots n,$ we would have $T_{g_i} A_i \in \bP'$. This is a contradiction, and completes the proof that $c_2$ is injective.

Next we show surjectivity of $c_2$, or, equivalently, surjectivity of $c$. Let $\bQ \in \GSpc \bK$; we will produce a prime $\bP \in \Spc \bK$ such that $c(\bP)=\bQ$. To do this, we will inductively define a collection of objects $\{A_1,..., A_i\}$ satisfying the following properties for $i \geq 2$:
\begin{enumerate}
    \item $A_1 \otimes A_2 \otimes ... \otimes A_i \not \in \bQ;$
    \item there exists a prime ideal $\bP_i \in \Spc \bK$ with $\bQ \subseteq \bP_i;$
    \item $A_j \not \in \bP_i$ for all $j=1,..., i$;
    \item $A_i \in \bP_{i-1}.$
\end{enumerate}
For the first step of the induction, one may simply pick any object $A_1$ which is not in $\bQ$; then since $\bQ$ is semiprime, we can pick a prime ideal $\bP_1$ containing $\bQ$ and not containing $A_1$.

Now we describe the inductive step (if it is possible). Suppose we have chosen objects $A_1,..., A_{i-1}$, such that $A_1\otimes... \otimes A_{i-1} \not \in \bQ$, and such that we have a prime ideal $\bP_{i-1} \in \Spc \bK$ with $\bQ \subseteq \bP_{i-1}$ and $A_j \not \in \bP_{i-1}$ for all $j=1,..., i-1$. Suppose there exists $A_i' \not \in \bQ$ satisfying $T_g A_i' \in \bP_{i-1}$ for all $g \in G$, that is, $A_i' \in c(\bP_{i-1}).$ Using \leref{gprime-char}, since both $A_1 \otimes... \otimes A_{i-1}$ and $A_i'$ are not in $\bQ$, there exists an object $B \in \bK$ and an element $g \in G$ such that
\[
A_1 \otimes A_2 \otimes... \otimes A_{i-1} \otimes B \otimes T_g A_i' \not \in \bQ.
\]
Set $A_i:=B \otimes T_g A_i'$. Note that $A_i \in \bP_{i-1}$, since $T_g A_i' \in \bP_{i-1}$, and that $A_1 \otimes ... \otimes A_i \not \in \bQ$ by construction. Then using again semiprimeness of $\bQ$, we can find a prime ideal $\bP_i$ with $\bQ \subseteq \bP_i$, and 
\[
A_1 \otimes A_2... \otimes A_i \not \in \bP_i;
\]
this implies that each $A_j \not \in \bP_i$, for $j=1,..., i.$ This completes the inductive step.

To complete the proof, we now consider the descending chain of closed subsets
\[
V(A_1) \supsetneq V(A_1) \cap V(A_2) \supsetneq ... \supsetneq V(A_1) \cap V(A_2) \cap... \cap V(A_i) \supsetneq...
\]
These inclusions are strict, since by construction we have $\bP_i \in V(A_1) \cap... V(A_i)$, but $A_{i+1} \in \bP_i$, that is, $\bP_i \not \in V(A_{i+1}).$ By Noetherianity, this process must terminate; that is, at some step, say $n$, we cannot find $A_{n+1}'$ with $A_{n+1}' \not \in \bQ$ and $T_g A_{n+1}' \in \bP_n$ for all $g \in G$. In other words, we cannot find $A_{n+1}' \in c(\bP_n)$ with $A_{n+1}' \not \in \bQ$, so $c(\bP_n) \subseteq \bQ$. But since $\bQ \subseteq \bP_n$ by construction, we have $c(\bP_n)=\bQ$. 
\end{proof}

We are now ready to complete the proof of the main theorem (\thref{maintheorem}). That is, by the previous result, we have a bijection between points of $\GSpc \bK$ and $G$-orbits of points of $\Spc \bK$; it remains to describe $\GSpc \bK$ as a topological space. We proceed to describe the closed sets of $\GSpc \bK$ in terms of orbits, and show that the specialization closed sets parametrize the $G$-ideals of $\bK$.

\bpr{topol-gspc}
Let $G$ be a group acting on a half-rigid $\Spc$-Noetherian monoidal triangulated category $\bK$. Using the bijection $c_2: (\Spc \bK)/G \to \GSpc \bK$ to identify points of $\GSpc \bK$ with orbits of points of $\Spc \bK$, if $W$ is a closed subset of $\Spc \bK$, then the orbit of $W$ in $(\Spc \bK)/G$ corresponds to a closed subset of $\GSpc \bK$. Furthermore, these closed subsets generate the topology on $\GSpc \bK$.   \epr

\begin{proof}
We check that
\begin{align*}
    V^G(A)&=\{ \bQ \in \GSpc \bK : A \not \in \bQ\}\\
    &=\{ c(\bP) : \bP \in \Spc \bK, A \not \in c(\bP)\}\\
    &= \{ c(\bP): \bP \in \Spc \bK, T_g A \not \in \bP \text{ for some } g \in G\}\\
    &=\{ c(\bP) : \bP \in \Spc \bK, A \not \in \bP\}\\
    &=\{ c(\bP) : \bP \in V(A)\}.
\end{align*}
The second equality follows since the core map $c$ is surjective, by \prref{bijective-core}. The fourth equality follows since $c(\bP)=c(T_g \bP)$ for any $g \in G$. Under the bijection $c_2$ between the sets $(\Spc \bK)/G$ and $\GSpc \bK$, we see that closed subsets of $\GSpc \bK$ are generated by the sets corresponding to the subsets 
\[
\{G.\bP : \bP \in V(A)\} \subseteq (\Spc \bK)/G.
\]
Recall that a general closed set $V^G(\bS)$, for a collection $\bS$ of objects of $\bK$, is the intersection of closed sets $V^G(A)$ for $A \in \bS$. 
\end{proof}

Note in particular that if $G$ is finite, or more generally if the ideal $\langle T_g A : g \in G\rangle$ is finitely-generated for all $A \in \bK$, then $\GSpc \bK$ is precisely the topological quotient of $\Spc \bK$ by the action of $G$. This recovers the results \cite[Theorem 9.2.3]{NVY1}, \cite[Section 10.C]{NP}, and  \cite[Theorem 9.1.1]{NVY3}, by \exref{group-hopf} and \prref{up-down-ideals}.

Lastly, we show that under our standard assumptions the $G$-ideals of $\bK$ are parametrized by specialization-closed subsets of $\GSpc \bK$. 

\bpr{ideal-class-gspc}
Let $G$ be a group acting on a half-rigid $\Spc$-Noetherian monoidal triangulated category $\bK$. The $G$-Balmer support $\Phi_G$ on $\GSpc \bK$ induces a bijection between $G$-ideals of $\bK$ and specialization-closed subsets of $\GSpc \bK$.
\epr

\begin{proof}
On one hand, by definition, if $\bI$ is a $G$-ideal of $\bK$, then its support is
\[
\Phi_G(\bI)= \bigcup_{A \in \bI} V^G(A),
\]
which is a specialization-closed set. The map which sends a specialization-closed subset $S$ of $\GSpc \bK$ to the $G$-ideal
\[
\bigcap_{\bQ \in \GSpc \bK \text{ with } \bQ \not \in S} \bQ
\]
is a left inverse to $\Phi_G$, using the fact that every $G$-ideal is an intersection of $G$-primes  by \coref{duals}. Indeed, it follows from \leref{g-ideal-supp} that $G$-ideals are in bijection with arbitrary unions of closed sets of the form $V^G(A)$, for $A \in \bK$.

It remains to be shown that every specialization-closed subset of $\GSpc \bK$ is a union of closed subsets of the form $V^G(A)$, that is, that an arbitrary closed subset $V^G(\bS)$, for $\bS$ a collection of objects of $\bK$, is a union of closed sets of the form $V^G(A)$. We set the notation that if $\bS$ is a collection of objects, then $G^{\bS}$ is the collection of set maps $\bS \to G$. For $\phi \in G^{\bS}$, we set $T_\phi(\bS)$ to be the collection of objects
\[
T_\phi(\bS):=\{ T_{\phi(A)} A : A \in \bS\}.
\]

 Note that by \thref{balmer-class}, we know that for any collection of objects $\bS'$, there exists some object, say $A_{\bS'}$, such that in $\Spc \bK$ we have $V(\bS') = V(A_{\bS'})$. Hence
\begin{align*}
    V^G(\bS) &=\bigcap_{A \in \bS} V^G(A)\\
    &= \bigcap_{A \in \bS} \{ c(\bP) : \bP \in V(A)\}\\
    &= \{ c(\bP) : \bP \in \Spc \bK, \text{ and} \; \;\forall A \in \bS, \; \; \exists g \in G \text{ with } T_g A \not \in \bP\}\\
    &=\bigcup_{\phi \in  G^{\bS}} \{ c(\bP): \bP \in V(T_\phi(\bS))\}\\
    &= \bigcup_{\phi \in G^{\bS}}\{ c(\bP): \bP \in V(A_{T_\phi(\bS)})\}\\ 
        &= \bigcup_{\phi \in G^{\bS}}V^G(A_{T_\phi(\bS)}).\\ 
\end{align*}
\end{proof}

Note that thick ideals are parametrized by specialization-closed sets and not by Thomason-closed subsets, since there is no guarantee that $\GSpc \bK$ is Noetherian. In particular, this means that if $\GSpc \bK$ is non-Noetherian (we provide an example below in \exref{perf-nonnoeth}), the $G$-ideal lattice of $\bK$ (equivalently, the ideal lattice of $\KrG$) cannot be equivalent to the ideal lattice of any symmetric or braided monoidal triangulated category, by \thref{balmer-class}. This gives potentially one of the first ``purely noncommutative" examples of Balmer spectra, in light of \cite[Question 4.2]{NP}.

\bex{nvy-gspec}
Let $\bC$ be a finite tensor category and $\bK:=\ul{\bC}$ its stable category. Recall that there is a comparison map $\rho: \Spc \bK \to \Proj C^\bullet_{\bK}$ (see \cite[Section 6.2]{NVY3}), where $C^\bullet_{\bK}$ is a particular subalgebra of the cohomology ring $\bigoplus_{i \geq 0} \Hom_{\bK} (\unit, \Sigma^i \unit)$ (see \cite[(1.1)]{NVY3}), under a weak finite-generation assumption (see \cite[Section 1.5]{NVY3}). Suppose that $\rho$ is an isomorphism; recall that this is conjectured to always hold \cite[Conjecture E]{NVY3}. Suppose also that $G$ is a finite group acting on $\bK$. By a similar argument to \cite[Theorem 5.2.1]{NVY3}, we have 
$C^\bullet_{\bK \rtimes G} \cong (C^\bullet_{\bK})^G.$ Then by \prref{up-down-ideals} and \prref{topol-gspc}, we have 
\[
\Spc (\bK \rtimes G) \cong \GSpc \bK \cong (\Proj C^\bullet_{\bK})/G \cong \Proj (C^\bullet_{\bK})^G \cong \Proj C^\bullet_{\bK \rtimes G}.
\]
That is, under our assumption that $\bC$ satisfies \cite[Conjecture E]{NVY3}, it follows that $\bC \rtimes G$ also satisfies \cite[Conjecture E]{NVY3}. This generalizes \cite[Theorem 9.1.1]{NVY3}, which proved the corresponding theorem for smash coproduct Hopf algebras. 
\eex

\bex{smallquant-gspec}
Let $\mf{g}$ be the Lie algebra of a simple algebraic group $G$ over an algebraically closed field $\kk$ of characteristic 0 and $u_{\zeta}(\mf{g})$ be the small quantum group at a root of unity, as in \exref{cocleft-action}. It is a conjecture of Negron and Pevtsova that $\Spc \ul{\modd}(u_{\zeta}(\mf{g})) = \mathbb{P}(\mc{N})$, where $\mc{N}$ is the nilpotent cone of $\mf{g}$ \cite[Section 17.4]{NP1}. This is proven for $G=\SL_2$ \cite[Corollary 17.2]{NP1}. In the $\SL_2$ case, $\mc{N} \cong \mathbb{A}_{\kk}^2$, so that $\Spc \modd(u_{\zeta}(\mf{g})) \cong \mathbb{P}^1_{\kk}$, and the action of $G$ is transitive on closed points. By \prref{up-down-ideals} and \prref{topol-gspc}, it follows that $\GSpc \ul{\modd}(u_{\zeta}(\mf{g})) \cong \Spc (\ul{\modd}(u_{\zeta}(\mf{g}) ) \rtimes G)$ has one closed point and one generic point. More generally, assuming the Negron--Pevtsova conjecture, closed points of $\GSpc \ul{\modd}(u_{\zeta}(\mf{g}))$ would be in bijection with nilpotent orbits, a finite set \cite[Section 2.5]{Jantzen}. 
\eex

\bex{perf-nonnoeth}
Let $X:=\mathbb{A}^2_{\mathbb{C}}=\Spec \mathbb{C}[x,y]$ be the affine scheme of complex affine 2-space. Recall that we have the rigid monoidal triangulated category $\Perf(X)$ as in \exref{perf}. Consider $G=\mathbb{Z}$ acting on $X$ via
\begin{align*}
    \mathbb{Z} & \rightarrow \Aut(\mathbb{C}[x,y])\\
    1 & \mapsto (x \mapsto x+1, y \mapsto y).
\end{align*}
Hence the action of $n \in \mathbb{Z}$ on a closed point $(a,b)=\langle x-a, y-b\rangle$ for $a,b \in \mathbb{C}$ is given by
\[
n: (a,b)\mapsto (a+n, b).
\]
We have an action of $\mathbb{Z}$ on $\Perf(X)$ by \exref{perf-action}, and the induced action of $\mathbb{Z}$ on $\Spc \Perf(X)$ corresponds to the original action of $G$ on $X$ when we identify $\Spc \Perf(X)$ and $X$ via the Thomason isomorphism. Let $V$ be the closed subset of $\mathbb{Z}$-$\Spc \Perf(X)$ corresponding to the orbit of the closed set $Z(x)$ of $\Spc \Perf(X)$, and $W$ be the closed subset of $\mathbb{Z}$-$\Spc \Perf(X)$ corresponding to the orbit of the closed set $Z(mx-y)$ of $\Spc \Perf(X)$, for some integer $m$. Then $W \cap V$ corresponds to the orbits of the closed points $(0, ma)$ for any $a \in \mathbb{Z}$; denote this closed subset $Y(m)$. Then we can construct an infinite descending chain
\[
Y(m) \supsetneq Y(m^2) \supsetneq Y(m^3) \supsetneq ...
\]
Hence $\mathbb{Z}$-$\Spc \Perf(X)$ is non-Noetherian, and so $\Perf(X) \rtimes \mathbb{Z}$ gives one of the first known examples of a monoidal triangulated category whose thick ideals are not parametrized by Thomason-closed subsets of its Balmer spectrum.
\eex



\end{document}